\documentclass{article}
\usepackage{ifthen}
\newboolean{isPreprint}
\setboolean{isPreprint}{true}
\newenvironment{keyword}{{\bf Keywords:}}{\vspace{1em}}

\newenvironment{frontmatter}{}{}
\date{}
\usepackage[dvipsnames,usenames]{xcolor}
\usepackage{amsmath}
\usepackage{amsfonts} 
\usepackage{latexsym}
\usepackage{algpseudocode}
\usepackage{algorithm}
\usepackage{listings}
\usepackage{graphicx}
\newcommand{\R}{\mathbb R}

\newcommand{\vect}[1]{\mathbf{#1}}
\newcommand{\tensor}[1]{\mathsf{#1}}

\newcommand{\diam}{\text{diam }}

\ifthenelse{\boolean{isPreprint}}
{\usepackage{amsthm}
\newtheorem{definition}{Definition}
}
{}

\begin{document}
\begin{frontmatter}
\title{A Massively Parallel Algebraic Multigrid Preconditioner based
  on Aggregation for
  Elliptic Problems with Heterogeneous Coefficients}

\ifthenelse{\boolean{isPreprint}}
{
\author{Markus Blatt\footnotemark[2]\ \footnotemark[1] \and Olaf Ippisch\footnotemark[2]
  \and Peter Bastian\footnotemark[2]}
\renewcommand{\thefootnote}{\fnsymbol{footnote}}
\footnotetext[1]{Dr. Markus Blatt - HPC-Simulation-Software \&
  Services, Hans-Bunte-Str. 8-10, D-69123 Heidelberg, Germany, email:
  markus@dr-blatt.de, URL: http://www.dr-blatt.de}
\footnotetext[2]{Interdisziplin\"ares  Zentrum f\"ur Wissenschaftliches Rechnen,
Ruprechts-Karls-Universit\"at Heidelberg, Im Neuenheimer Feld 368, D-69120 Heidelberg,
 Germany}
\maketitle
}
{\author[mb]{Markus Blatt\corref{cor1}}
\ead{markus@dr-blatt.de}
\ead[url]{http://www.dr-blatt.de/}
\cortext[cor1]{Corresponding author}
\author[hd]{Olaf Ippisch}
\author[hd]{Peter Bastian}

\address[mb]{Dr. Markus Blatt - HPC-Simulation-Software \&
  Services, Hans-Bunte-Str. 8-10, D-69123 Heidelberg, Germany}
\address[hd]{Interdisziplin\"ares  Zentrum f\"ur Wissenschaftliches Rechnen,
Ruprechts-Karls-Universit\"at Heidelberg, Im Neuenheimer Feld 368,
D-69120 Heidelberg, Germany}
}
%\listoftodos
\begin{abstract}
This paper describes a massively parallel algebraic multigrid method based
on non-smoothed aggregation. 
A greedy algorithm for the aggregation combined with an appropriate 
strength-of-connection criterion makes it especially suited for solving
heterogeneous elliptic problems. 
Using decoupled aggregation on
each process with data agglomeration onto fewer processes on the
coarse level, it weakly scales well in terms of both total time to
solution and time per iteration to nearly 300,000 cores. Because of
simple piecewise constant interpolation between the levels, its memory
consumption is low and allows solving problems with more than
$10^{11}$ degrees of freedom.
\end{abstract}
%\begin{keywords} 
%\begin{AMS}
%  65F08, 65N08, 65N55, 65Y05
%\end{AMS}
\end{frontmatter}
\pagestyle{myheadings}
\thispagestyle{plain}
\markboth{Blatt, Ippisch, Bastian}{MASSIVELY PARALLEL AGGREGATION AMG}{}

\section{Introduction}
\label{sec:introduction}

When solving elliptic or parabolic partial differential equations (PDEs)
most of the 
computation time is often spent in solving the arising linear
algebraic equations. This
demands for highly scalable parallel solvers capable of running on
recent supercomputer.  The current
trend in the development of high performance supercomputers is to
build machines that utilize  more and more cores with less memory
per core, but interconnected with low latency networks. To be able to
still solve problems of reasonable size the parallel linear solvers
need to be (weakly) scalable and have a very small memory footprint.

Besides domain decomposition methods the most scalable and fastest methods are multigrid
methods. They can solve these
linear systems with optimal or nearly optimal complexity, i.e. at most $O(N\log N)$ operations for
$N$ unknowns. 
Among them algebraic multigrid methods (AMG) are
especially suited for problems with heterogeneous or anisotropic
coefficient tensors on
unstructured grids. They build a hierarchy of matrices using their
graphs and thus adapt the coarsening to the problem solved.

Parallel geometric multigrid implementations exist since at least 25
years, see e.g. \cite{Frederickson:1987:PSM:867287}. 
Since about 15 years, several parallel algebraic multigrid codes have been developed
\cite{BoomerAMG:1999,Tuminaro00parallelsmoothed,Krechel20011009,Fudap_algebraic:2005,Joubert:2006}.
%The most recent developments include also parallel algebraic multigrid on GPUs
%\cite{Haase:2009:PAM:2127669.2127675,Gundolf_AMG_GPU_latest:2010}.
Classical AMG \cite{ruge87:multig_method_amg}
divides the fine level unknowns into two groups: the ones also
represented on the coarse level, and the ones that exist only on the
fine level.  
In parallel versions of the coarsening algorithm neighboring processors
have to agree on the coarse/fine splitting in the overlap region
without spending too much time in communication while still ensuring
small work per cycle and good convergence rate. Alber and Olson \cite{alber07:_paral_coars_grid_selec}
give a comprehensive comparison of several parallel coarsening algorithms
that indicate that achieving small work per cycle and good convergence rate
is difficult to achieve at the same time.
Adapted coarsening heuristics, aggressive coarsening strategies as well
as hybrid implementations (shared memory on a node, message passing between nodes) have been
developed to overcome this problem
\cite{Stueben_1999,Sterck:2006,alber07:_paral_coars_grid_selec,Baker:2011:CSA:2058524.2059558}.
%\cite{alber07:_paral_coars_grid_selec,Griebel.Metsch.Schweitzer:2007,sterck08:_distan_two_inter_for_paral_algeb_multig}

AMG based on aggregation, see \cite{Braess,Vanek,Raw_1996},
 clusters the fine level
unknowns into aggregates. Each aggregate represents an unknown on the
coarse level and its basis function is a linear combination of the
fine level basis 
functions associated with the aggregate. Two main classes of the
method exist. Non-smoothed aggregation AMG, see
\cite{Braess,Raw_1996,notay10,blattamg}, which uses simple piecewise 
constant interpolation, and smoothed aggregation AMG, that increases
interpolation accuracy by smoothing the tentative piecewise constant
interpolation. For the parallel versions of both classes no growth in
operator complexity is observed for increasing numbers of processes
\cite{Tuminaro00parallelsmoothed,notay10,blattamg}. Still the
smoothing of the interpolation operators increases the stencil size of
the coarse level matrices compared to the non-smoothed
version. Moreover, the non-smoothed version can be used in  straight
forward way for many systems of PDEs, see \cite{blattamg}.

In
this paper we describe a parallel AMG method  that uses a greedy
heuristic algorithm for the aggregation based on a strength of
connection criterion. This allows for building
round aggregates of nearly arbitrary size that do not cross high
contrast coefficient jumps. We use simple piecewise constant
interpolation between the levels to prevent an increase of the size of the
coarse level stencils.  Together with an implementation of the parallel linear
algebra operations based on index sets this makes the algorithm very scalable regarding the time
needed per iteration. Even though the number of iterations needed for
convergence does increase during weak scalability tests, the time
to solution is still very
scalable. We present numerical evidence that the approach is scalable up to
262,144 cores for realistic problems with highly variable coefficients.
At the same time the memory requirement of the algorithm is far less
than that of classical AMG methods.  This allows us to solve problems with more than
$10^{11}$ degrees of unknowns on an IBM Blue Gene/P using 64 racks.

We will start the paper in the next section with a description of the
algebraic multigrid method together with our heuristic greedy
aggregation algorithm for coarsening the linear systems. In Section
\ref{sec:parallelization} we describe the parallelization of the
algebraic multigrid solver and its components, namely the data
decomposition, smoothers, interpolation operators, and linear
operators. After presenting implementational details about the
parallelization and linear algebra data structures in Section
\ref{sec:impl-deta}, we conduct scalability tests of our method on
an IBM Blue Gene/P and an off-the-shelf multicore Linux cluster in Section
\ref{sec:numerical-results}. Our summary and conclusions can be found
in Section \ref{sec:summary}.

\section{Algebraic Multigrid}
\label{sec:algebr-mult-based}

The notation of parallel linear algebra algorithms can be simplified considerably
by the use of non-consecutive index sets. 
This allows one to use a single, global index set without the need for remapping indices
to a consecutive index set for each processor.
Following \cite{Hackbusch:1994},
for any 
finite index set $I\subset\mathbb{N}$ we define the vector space 
$\mathbb{R}^I$ to be isomorphic to $\mathbb{R}^{|I|}$ with components
indexed by $i\in I$. Thus $\vect x \in \mathbb{R}^I$ can be interpreted as 
a mapping $\vect x : I \to \mathbb{R}$ and $(\vect x)_i = \vect x(i)$. In the same 
way, for any two finite index sets $I, J\subset\mathbb{N}$ we write 
$\tensor A \in \mathbb{R}^{I\times J}$ with the interpretation 
$\tensor A : I \times J \to \mathbb{R}$ and $(\tensor A)_{i,j} = \tensor A(i,j)$. 
Finally, for any subset $I'\subseteq I$ we define the restriction
matrix
\begin{equation}
\label{eq:restrict} 
\tensor R_{I,I'} : \mathbb{R}^I \to \mathbb{R}^{I'} \text{\quad as \quad} 
 (\tensor R_{I,I'} \vect x)_i = (\vect x)_i \ \forall i\in I'
\end{equation}
(which
  corresponds to simple injection). 

On a given domain $\Omega$ we are interested in solving the model
problem 
\begin{equation}
  \label{eq:model}
  \nabla \cdot (\tensor K(x) \nabla u) = f\,,\quad \text{on } \Omega
\end{equation}
together with appropriate boundary conditions. Here, the symmetric
positive definite tensor $\tensor K(x)$, dependent on the position $x$
within the domain $\Omega$, is allowed to be discontinuous. Given an
admissible mesh $T_h$ that for simplicity resolves the boundary and
possible discontinuities in the tensor $\tensor K(x)$, discretizing
\eqref{eq:model} using conforming lowest order Galerkin finite element or finite volume
methods yields a linear system
\begin{equation}
  \label{eq:ls}
  \tensor A x = b\,,
\end{equation}
where $\tensor A: \R^I \to \R^I$ is the linear operator, and $\vect x,
\vect b\in\R^I$
are vectors. For an extension to discontinuous Galerkin methods see
\cite{blatt12_discGalerkin}. We strive to solve this linear system using our algebraic
multigrid method described below.

The excellent computational complexity of multgrid methods is due to the
following main idea. Applying a few steps of a smoother (such as
Jacobi or Gauss-Seidel) 
to the linear system usually leads to a smooth
error that cannot be reduced well by further smoothing. Given a
prolongation operator $\tensor P$ from a coarser linear
system, this error is  then further reduced 
using a correction $u^\text{coarse}$ on a coarser linear system $\tensor P^T\tensor A\tensor P \vect u^\text{coarse}
= \tensor P^T (\vect b-\tensor A \vect x)$. We use the  heuristic
algorithm presented in Subsection 
\ref{sec:coars-algor-based} to build the prolongation operator
$\tensor P$. If the system is already small enough, we solve it using a
direct solver. Otherwise we recursively apply a few steps of the smoother
and proceed to an even coarser linear system until the size of the coarsest level
is suitable for a direct solver. After applying the coarse level
solver, we prolongate the correction to the next finer level, add it to the
current guess, and apply a few steps of the smoother.

\subsection{Coarsening by Aggregation}
\label{sec:coars-algor-based}

To define the prolongation operator $\tensor P$ we rely on a greedy
and heuristic aggregation algorithm, that uses the graph of the matrix
as input. It is an extension of the 
version published by Raw (cf. \cite{raw85:amg}) for algebraic
multigrid methods (see  also \cite{scheichl}). 

Let $G=({\mathcal V},{\mathcal E})$ be a graph with a set of vertices
${\mathcal V}$ and edges ${\mathcal E}$
and let  $w_{\mathcal E} \,:\, {\mathcal E} \rightarrow \R$ and
$w_{\mathcal V}\,:\, {\mathcal V} \rightarrow \R$ be positive weight
functions. For the examples in this paper $w_{\mathcal
  E}((i,j))=\frac{1}{2}((\tensor A)_{ij}-|\,(\tensor A)_{ij}|)$ is used, i.e. $0$ for positive
off-diagonals and the absolute value otherwise,
and $w_{\mathcal V}(i)=(\tensor A)_{ii}$. For matrices arising from the discretisation of systems of PDE, for which our aggregation scheme is applicable as well,
$w_{\mathcal E}$ and $w_{\mathcal V}$ could e.g. be the row-sum norm
of a matrix block (see e.g. \cite{blattamg}). These functions are used
to classify the edges and vertices 
of 
our graph. Let 
$$N(i):=\{j \in {\mathcal V}\mid \exists (j,i) \in {\mathcal E}\}$$
 be the set of adjacent vertices of vertex $i$,
let
$$N(a):=\{j \in {\mathcal V}\setminus a\mid \exists k \in a
\text{with} i \in N(k) \}$$
be the set of adjacent vertices of a set of vertices,  and let
 \begin{equation}
    \label{eq:symmetric_strong1}
    \eta_{\max}(i) := \max_{k \in N(i)} 
\frac{w_{\mathcal E}((k,i))\: w_{\mathcal E}((i,k))}{w_{\mathcal V}(i)\:w_{\mathcal V}(k)}.
%&\quad\mbox{for symmetric } \tensor{A}\\
 %   \label{eq:strong1}
  %  \eta_{\max}(i):=& \max_{k \in N(i)} \frac{-w_E((k, i))}{w_V(i)}&\quad\mbox{otherwise}\,,
  \end{equation}

{
\renewcommand{\labelenumi}{(\alph{enumi})}
 \begin{enumerate}
\item  An edge $(j, i)$ is called {\em strong}\index{matrix~graph!strong
  edge}, if and only if
  \begin{equation}
    \label{eq:symmetric_strong2}
    \frac{w_{\mathcal E}((i, j))\:w_{\mathcal E}((j, i))} {w_{\mathcal V}(i)\: w_{\mathcal V}(j)} > \delta
    \:\min(\eta_{\max}(i), \eta_{\max}(j)),%\quad&\mbox{for symmetric }\tensor{A}\\
    %\label{eq:strong2}
    %\frac{-w_E((j, i))}{w_V(i)} > \alpha \:\eta_{\max}(i)\quad&\mbox{otherwise}
  \end{equation}
  for a given threshold $0<\delta<1$. We denote by $N_{\delta}(i)\subset N(i)$ 
  the set of all vertices adjacent to $i$ that are connected to it via a 
  strong edge. Furthermore we call an edge $(i,j)$ a one-way strong
  connection if edge $(j,i)$ is not strong. If both $(i,j)$ and
  $(j,i)$ are strong we call the edges a two-way strong connection.

\item A vertex $i$ is called {\em
  isolated}\index{matrix~graph!isolated vertex} if and only if $\eta_{\max}(i) <
\beta$, for a prescribed
  threshold $0 < \beta \ll 1$. We denote by $\mathit{ISO}(\mathcal V)\subset \mathcal V$ the set of all 
  isolated vertices of the graph.
\end{enumerate}
}

{For symmetric positive definite M-matrices arising from
  problems with constant diffusion coefficients our strength of
  connection criterion is similar to the traditional ones for the AMG
  of Ruge and St\"uben \cite{ruge87:multig_method_amg} and for
  smoothed aggregation \cite{PVanek_JMandel_MBrezina_1996a}. For
  non-symmetric matrices or problems with discontinuous coefficients 
  it differs from  them. It is
  especially tailored for the latter. At the interfaces of the jumps the Ruge St\"uben
  criterion might classify a connection between two vertices strong
  in one direction and weak in the other one. Actually, no aggregation
  should happen across this interface. The smoothed aggregation
  criterion falsely classifies positive off-diagonal values as strong
  while ours does not do this with an appropriate weight function. For
more details see \cite{blattamg}.}

Our greedy aggregation algorithm is described in Algorithm
\ref{alg:sequential_coarsen}. 
\begin{algorithm}[H]
  \caption{Build  Aggregates}\label{alg:sequential_coarsen}
  \begin{algorithmic}[0]
    \Procedure{Aggregation}{${\mathcal V}$, ${\mathcal E}$, $s_\text{min}$, $s_\text{max}$,
      $d_\text{max}$} %\Comment{Arguments: graph, minimum and maximum aggregate size, maximum radius}
    \State $U \gets \{ v \in {\mathcal V} \setminus
    \mathit{ISO}({\mathcal V}) : \text{v not on Dirichlet boundary}\}$
    \Comment{First candidates are non-isolated vertices}
    \State $I \gets \emptyset$ \Comment{Coarse index set}
     \State $S \gets \{u \in U : |N_\text{na}(u)| \leq
    |N_\text{na}(w)| \ \forall w \in U\}$  %\Comment{Seed stack}
    \State Select $v \in S$
    \While{$U\not=\emptyset$}
      \State $a_v \gets \{v\}$ \Comment{Initialize new aggregate}
      \State $U \gets U \setminus a_v$
      \State $I \gets I \cup \{v\}$
      \State \Call{growAggregate}{$a_v$, ${\mathcal V}$, ${\mathcal E}$, $s_\text{min}$, $d_\text{max}$, $U$}
      \State \Call{roundAggregate}{$a_v$, ${\mathcal V}$, ${\mathcal E}$, $s_\text{max}$, $U$}
      \If{$|a_v|=1$} 
      \Comment{Merge one vertex aggregate with neighbor}
        \State $C \gets \{a_j : j\in I\setminus\{v\} \ \text{and} \
        \exists w \in a_j
        \text{ with } w  \in N_\delta(v)\}$
        \If{$C\not=\emptyset$}
        \State Choose $a_k\in C$
        \State $I \gets I \setminus \{v\}$
        \State $a_k \gets a_k \cup a_v$
        \EndIf
      \EndIf
      \State $S \gets \{w: w \in N(a_v)\}$
      \If{$U\not=\emptyset$}
        \If{$S=\emptyset$}
        \State $S \gets \{u \in U : |N_\text{na}(u)| \leq
    |N_\text{na}(w)| \ \forall w \in U\}$
        \EndIf
        \State Select $v \in S$
      \EndIf
    \EndWhile
    \State $U \gets \{v \in \mathit{ISO}({\mathcal V}) : \text{v not on Dirichlet boundary}\}$
    \Comment{Aggregate isolated vertices}
    \While{$U\not=\emptyset$}
    \State Select arbitrary seed $v \in U$
    \State $a_v \gets \{v\}$
    \State $U \gets U \setminus a_v$
    \State $I \gets I \cup \{v\}$
    \State \Call{GrowIsoAggregate}{$a_v$, ${\mathcal V}$, ${\mathcal E}$,
      $s_\text{min}$, $d_\text{max}$, $U$}
    \EndWhile
    \State $\mathcal A \gets \{a_i : i \in I\}$
    \State \Return $(\mathcal A, I)$
    \EndProcedure
  \end{algorithmic}
\end{algorithm}
Until all non-isolated vertices, that are not part of a Dirichlet boundary, are
aggregated, we start a new aggregate $a_v$ with a non-isolated vertex
$v$. We prefer vertices as seeds which have
the least number of non-aggregated neighbors $N_\text{na}$. The index
of the seed vertex is associated 
with this new aggregate and added to the index set $I$. The algorithm
returns both the index 
set $I$ for the set of aggregates as well as the set $\mathcal A = 
\{a_i : i \in I\}$ of all aggregates it has built.

The first step in the construction of an aggregate in Algorithm
\ref{alg:sequential_coarsen} is to add new vertices to our aggregate until we 
reach the minimal prescribed aggregate size $s_\text{min}$. This is 
outlined in Algorithm \ref{alg:grow_aggre}. Here $\diam(a,v)$ denotes
the graph diameter of the subgraph
$G=({\widetilde{\mathcal V}},{\widetilde{\mathcal E}})$
with ${\widetilde{\mathcal V}}=a\cup\{v\}$ and $\widetilde{\mathcal
    E} =\{(i,j) \in {\mathcal E} \mid i,j \in {\widetilde{\mathcal
    V}}\}$. Recall that the diameter of a graph is the longest
shortest path between any two vertices of the graph.
\begin{algorithm}[H] \caption{Grow Aggregate Step}\label{alg:grow_aggre}
  \begin{algorithmic}[0]
    \Function{GrowAggregate}{$a$, ${\mathcal V}$, ${\mathcal E}$, $s_\text{min}$, $d_\text{max}$, $U$}
    \While{$|a|\leq s_\text{min}$}\Comment{Makes
      aggregate $a$ bigger until its size is $s_\text{min}$}
       \State $C_0 \gets \{v \in N(a) : \text{ diam}(a,v)\leq d_\text{max}\}$ \Comment{Limit the diameter of the aggregate}
       \State $C_1 \gets \{v \in C_0 : \text{ cons}_2(v,a) \geq
       \text{ cons}_2(w,a) \ \forall w \in N(a)\}$
       \If{$C_1=\emptyset$}\Comment{No candidate with two-way connections}
         \State $C_1 \gets \{v \in C_0 : \text{ cons}_1(v,a) \geq
       \text{ cons}_1(w,a) \ \forall w \in N(a)\}$
       \EndIf
       \If{$|C_1| >1$}\Comment{More than one candidate}
         \State $C_1 \gets \{ v \in C_1 : \frac{\text{ connect}(v,
         a)}{|N(v)|} \geq \frac{\text{ connect}(w, a)}{|N(w)|} \ \forall w \in
       C_1\}$
       \EndIf
       \If{$|C_1| >1$}\Comment{More than one candidate}
         \State $C_1 \gets \{ v \in C_1 : \text{ neighbors}(v,
         a) \geq \text{ neighbors}(w, a) \ \forall w \in
       C_1\}$
       \EndIf
       \If{$C_1=\emptyset$}
          \textbf{break}
       \EndIf
       \State Select one candidate $c\in C_1$
       \State $a \gets a \cup \{c\}$ \Comment{Add
         candidate to aggregate}
       \State $U \gets U\setminus \{c\}$
     \EndWhile

   \EndFunction
  \end{algorithmic}
\end{algorithm}
When adding new vertices, we always
choose a vertex within the prescribed maximum diameter $d_\text{max}$ of the aggregate which has the highest number of strong connections to the vertices already in
the aggregate. Here we give preference to vertices where both edge $(i,j)$ 
and edge $(j,i)$ are strong. 
The functions $\text{cons}_1(v, a)$ and $\text{cons}_2(v,
a)$ return the number of one-way and two-way
strong connections between
the vertex $v$ and all vertices of the aggregate $a$,
respectively. Note the for the examples and the choice of $w_{\mathcal
  E}$ in this paper all strong edges are two-way strong connections.

If there is more than one candidate, we want to choose a vertex with a
high proportion of strong connections to other vertices not belonging
to the current aggregate, while favoring connections to vertices
which belong to aggregates that are already connected to the current
aggregate. We therefore define a function $\text{connect}(v,a)$, which
counts neighbors of $v$. Neighbors of $v$ that are not yet aggregated or belong to an 
aggregate that is not yet connected to aggregate $a$ are counted once.
Neighbors of $v$ that belong to aggregates that are already connected
to aggregate $a$ are counted twice. 

If there is still more than one candidate which maximizes
$\frac{\text{connect}(v,a)}{|N(v)|}$, we 
choose the candidate, which has the maximal 
number $\text{neighbors}(v,a)$ of neighbors of vertex $v$ that are not
yet aggregated 
neighbors of the aggregate $a$. This criterion tries to
maximize the number of candidates for choosing the next vertex.

\begin{algorithm} [H]
  \caption{Round Aggregate Step}\label{alg:round_aggre}
  \begin{algorithmic}[0]
    \Function{RoundAggregate}{$a$, ${\mathcal V}$, ${\mathcal
        E}$, $s_\text{max}$, $U$}
    \While{$|a|\leq s_\text{max}$}\Comment{Rounds aggregate $a$ while size $< s_\text{max}$}
      \State $D\gets \{w \in N(a)\cap U : \text{ cons}_2(w,
      a)>0 \text{ or cons}_1(w, a)>0\}$ 
      \State $C\gets \{v \in D : |\{N(v)\cap U\}| > |\{N(v) \cap U\}|\}$
      \State Select arbitrary candidate $c \in C$
      \State $a\gets a \cup \{c\}$ \Comment{Add
         candidate to aggregate}
       \State $U \gets U\setminus \{c\}$
    \EndWhile
    \EndFunction
  \end{algorithmic}
\end{algorithm}  

In a second step we aim to make the aggregates
``rounder''. This is sketched in Algorithm~\ref{alg:round_aggre}. 
We add all non-aggregated adjacent vertices that have more connections to the 
current aggregate than to other non-aggregated vertices until we reach the maximum allowed size $s_\text{max}$ of our aggregate. 

If after these two steps an aggregate still consists of only one
vertex, we try to find another aggregate that the vertex is strongly
connected to. If such an aggregate exists, we add the vertex to that
aggregate and choose a new seed vertex.

Finally, once all the non-isolated vertices are aggregated, we try to 
build aggregates for the isolated vertices, that are not part of a
Dirichlet boundary. Those kind of vertices can be produced on coarser
levels by successively aggregating in small regions bounded by
coefficient jumps. Where possible, we build 
aggregates of adjacent isolated vertices that have at least one 
common neighboring aggregate consisting of non-isolated vertices. 
This is done in the function {\sc \lstinline!GrowIsoAggregate!} which we do
not present here.  Our
aggregation algorithm ensures reasonable 
coarsening rates and operator complexities in this case.

Given the aggregate information $\mathcal A$, we define the piecewise
constant prolongation operator $\tensor P$ by
\begin{equation}
\label{eq:prolongation}
\tensor P(i,j)=\left\{
\begin{array}{rc}
  1 & \text{if } j \in a_i\\
  0 & \text{else}
\end{array}\right.
\end{equation}
and define the coarse level matrix using a Galerkin product as
$$
\tensor A^\text{coarse}= \frac{1}{\omega} \tensor P^T \tensor A\tensor
P\,.
$$
Note that the over-correctionfactor $\omega$  is
needed to improve the approximation properties of the coarse
correction. According to \cite{Braess} $\omega=1.6$ is a good default
and often sufficient for good convergence.

\section{Parallelization}
\label{sec:parallelization}

\subsection{Data Decomposition and Local Data Structures}
\label{sec:data-decomposition}

The most important and computationally expensive part in both
algebraic multigrid methods and many iterative solvers (especially
Krylov methods  and stationary iterative methods) is the application
of linear operators (in the simplest case realized by matrix multiplication). 
Therefore their construction is crucial. It
has to be made in such that it allows for the efficient application of the
linear operator \textit{and} the preconditioners based on these operators.

Let $I\subset\mathbb{N}$ be our finite index set and
let $\bigcup_{p\in\mathcal P} I_{(p)}$
be a (non-overlapping) partitioning of the index set for the
processes $\mathcal P$. This partitioning might be given by an external
partitioning software or by a parallel grid manager. 
In some cases these tools may not provide a partitioning but an overlapping
decomposition. In this case one can easily compute a partitioning
following \cite{ISTLParallel}.

In order to avoid communication during matrix multiplication
the index set $I_{(p)}$ is extended to
$$
{\widetilde I_{(p)}} = I_{(p)} \cup \{j \in I | \exists i \in I_{(p)}
\,:\,(\tensor A)_{ij}\not=0 \}\,,
$$
which are all indices of $x$ that are needed for computing the components
of the product $Ax$ corresponding to $I_{p}$ in process $p$.
Furthermore, we prefer to store quadratic matrices in each process in order
to represent the linear operator. To that end
each process $p$ stores $\widetilde{\tensor A}_{(p)} \in \R^{\widetilde
  I_{(p)}\times\widetilde I_{(p)}}$ and $\widetilde{ \vect x}_{(p)},
\widetilde{\vect b}_{(p)} \in \R^{\widetilde I_{(p)}}$. 
\begin{definition}
We call
$\widetilde{\vect x}_{(p)}$ stored {\em consistently} if $\widetilde{\vect x}_{(p)} = \tensor R_{I,\widetilde I_{(p)}} \vect x$ for all
$p\in\mathcal P$, with $\tensor R$ defined by \eqref{eq:restrict}. If  $\widetilde{\vect x}_{(p)} =\tensor R^T_{\widetilde
  I_{(p)},I_{(p)}} \tensor R_{I,I_{(p)}}\vect x$ holds
for all
$p\in\mathcal P$, we denote $\widetilde{\vect x}_{(p)}$ as {\em uniquely} stored.
\end{definition}

Let $\tensor A \in \R^{I \times I}$ be a global linear operator. Then on
process $p$ the local linear operator $\widetilde{\tensor A}_{(p)}$ stores the values
\begin{equation}
\label{eq:tilde_A_bar_tilde}
(\widetilde{\tensor A}_{(p)})_{ij} = \left\{
  \begin{array}{cl}
    (\tensor A)_{ij} & \text{if }
    i \in I_{(p)} \\
    \delta_{i,j} & \text{else}
  \end{array}
\right.
\end{equation}
where $\delta_{i,j}$ denotes the usual Kronecker delta. Denoting
$\tensor A_{(p)}=\tensor R_{I,I_{(p)}}\tensor A \tensor R^T_{I,I_{(p)}}$ and reordering the indices locally, such that
$i<j$ for all $i\in \widetilde I_{(p)}$ and $j\not\in \widetilde I_{(p)}
\setminus I_{(p)}$, $\widetilde{\tensor A}_{(p)}$ has the
following structure:
\begin{equation*}
  \begin{array}{ll}
    \widetilde{I}_{(p)}\left\{
      \begin{array}{l}
        I_{(p)}\left\{
          \begin{array}{l}
            \\\\\\
          \end{array}
        \right.
        \\
        \begin{array}{l}
        \end{array}
      \end{array}
    \right.
& \begin{array}{|ccc|c|}
  \hline
  &&&\\
  \text{     }&\tensor A_{(p)}&\text{     }&\;\ast\;\\
  &&&\\\hline
  &0&&I\\\hline
\end{array}
\end{array}\,.
\end{equation*}

Using this storage scheme for the local linear operator $\widetilde{\tensor A}_{(p)}$ and
applying it to a local vector $\widetilde{\vect x}_{(p)}$, stored consistently,
$(\widetilde{\tensor A}_{(p)} \widetilde{\vect x}_{(p)})(i)=(\tensor A\vect x)(i)$ holds for all $i\in
I_{(p)}$. Therefore the global application 
of the linear operator can be represented by computing
\begin{equation}\label{eq:parlinop}
\tensor A \vect x = \sum_{p\in\mathcal
  P}{\tensor R}^T_{I,I_{(p)}} {\tensor R}_{\widetilde
  I_{(p)},I_{(p)}}\left( \widetilde{\tensor A}_{(p)} {\tensor R}_{I,\widetilde I_{(p)}} \vect x\right)\,.
\end{equation}
Here the operators in front of the brackets represent a restriction of the results of the local computation to the (consistent) representation on $\R^{I_{(p)}}$ then a 
prolongation to global
representations and a summation of all these global representations. In contrast
to the notation used there is no
global summation and thus no global communication needed. It suffices
that every process adds only 
entries from other processes that actually store data associated with
indices in $\widetilde I_{(p)}$. Therefore this represents a next neighbor
communication followed by a local summation.

\subsection{Parallel Smoothers}
\label{sec:smoothers}

As smoothers we only consider so-called hybrid smoothers \cite{hybrid}. These
can be seen as block-Jacobi smoothers where the blocks are the
matrices $\tensor A_{(p)}$. Instead of directly solving the block systems a few
steps of a sequential smoother (e.g. Gauss-Seidel for hybrid
Gauss-Seidel) is applied. We always use only one step.

Let $\tensor M_{(p)}\in \R^{I_{(p)}\times I_{(p)}}$, $p \in \mathcal P$, be the sequential
smoother computed for matrix $\tensor A_{(p)}$, and
$\widetilde{\vect d}_{(p)}$ the defect. Then the consistently stored update
$\widetilde{\vect v}_{(p)}$ is computed by
applying the parallel preconditioner as
\begin{equation*}
  \widetilde{\vect v}_{(p)} = \tensor R_{I,\widetilde I_{(p)}}\sum_{p \in \mathcal P}
  \tensor R_{I,I_{(p)}}^T \left(\tensor M^k_{(p)} \tensor R_{\widetilde I_{(p)}, I_{(p)}}
    \widetilde d_{(p)}\right)\,,
\end{equation*}
for $k\geq 1$. This means that due to our storage according to
\eqref{eq:tilde_A_bar_tilde} we can apply multiple steps a local
preconditioner in a hybrid smoother without adding further communication.
Again as in the parallel linear operator \ref{eq:parlinop} the summation requires only a communication with processes, which share data associated with $\widetilde I_{(p)}$, and thus can be handled very efficiently.

\subsection{Parallel Coarsening}
\label{sec:parallel-aggregation}

The parallelization of the coarsening algorithm described in Section
\ref{sec:coars-algor-based} is
rather straightforward. It is simple and massively parallel
since the aggregation only occurs on vertices of the graph of
matrix $\tensor A_{(p)}$. Using this approach, the coarsening process will of course
deal better with the algebraic smoothness if the disjoint
matrix $\tensor A_{(p)}$ is split along weak edges. 
% This will not be possible in
% all cases and therefore agglomerating the data onto fewer processors will
% to some extend overcome this drawback.

The parallel approach is described in Algorithm
\ref{alg:parallel_coarsen}. It builds the aggregates
$\widetilde{\mathcal A}_{(p)}$ of this level and the
parallel index sets $\widetilde I_{(p)}^\text{coarse}$ for the next level in parallel. The parameters are the edges and vertices
of the matrix graph $G(\widetilde{\tensor A}_{(p)})=(
\widetilde I_{(p)}, \widetilde{\mathcal E}_{(p)})$ and the disjoint index set $I_{(p)}$. The rest of the
parameters are the same as for the sequential Algorithm \ref{alg:sequential_coarsen}. As a first step
a subset $(I_{(p)}, {\mathcal E}_{(p)})$ of the input graph that corresponds to the index
set $I_{(p)}$ is created. Then the
sequential aggregation algorithm is executed on this sub-graph. Based
on the outcome of this aggregation a map between indices and
corresponding aggregate indices is built and the information is
published to all other processes that share vertices of the
overlapping graph. Now every process knows the aggregate index of each
vertex of its part of the overlapping graph and constructs the
overlapping coarse index set and the aggregates. Note that this
algorithm only needs one communication step per level with the direct neighbors.

\begin{algorithm}
  \caption{Parallel Aggregation}
  \label{alg:parallel_coarsen}
  \begin{algorithmic}
    \Procedure{ParallelAggregation}{$I_{(p)}$, $\widetilde{\mathcal E}_{(p)}$,
      $s_\text{min}$, $s_\text{max}$, $d_\text{max}$}
    \State \textbf{On process $p\in\mathcal P$:}
%    \State ${\mathcal V}_{(p)} \gets \{v_k \in {\mathcal V}_{\widetilde I_{(p)}} \mid k \in I_{(p)}\}$ \Comment{Only
%      vertices owned by $p$}
    \State ${\mathcal E}_{(p)} \gets \{(k,l) \in \widetilde{\mathcal E}_{(p)} \mid k \in I_{(p)}\,,
    l \in I_{(p)}\}$ 
    \Comment{Only edges within $I_{(p)}$}
    %\State $E^O_{(p)} \gets  \widetilde E_{(p)} \setminus E_{(p)}$ \State $V^O_{(p)} \gets  \widetilde V_{(p)} \setminus V_{(p)}$
    \State $(I^\text{coarse}_{(p)}, \mathcal A_{(p)}) \gets$\Call{Aggregation}{$I_{(p)}$, ${\mathcal E}_{(p)}$, $s_\text{min}$, $s_\text{max}$,
      $d_\text{max}$}
    \State $\widetilde{\vect m}_{(p)}\gets 0 \in  \R^{\widetilde I_{(p)}}$
    \For{$a_k \in \mathcal A_{(p)}$}
    \State $(\tensor R^T_{I,\widetilde I_{(p)}}\widetilde{\vect m}_{(p)})_j \gets k \quad \forall
    v_j \in a_k$
    \EndFor
    \State $\widetilde{\vect m}_{(p)} \gets R_{I,\widetilde I_{(p)}} \sum_{q\in \mathcal P} \tensor
    R^T_{I, \widetilde I_{(q)}} \widetilde{\vect m}_{(q)}$ \Comment{Communicate aggregates mapping}
    \State $\widetilde I^\text{coarse}_{(p)} \gets \{k \mid
    \exists j \in \widetilde I_{(p)} \text{ with } (
    R^T_{I,\widetilde I_{(p)}}\widetilde{\vect m}_{(p)})_j
    =k\}$\Comment{Build coarse index set}
    \For{$k \in \widetilde I^\text{coarse}_{(p)}$}
    \State $a_k \gets\{ j \in \widetilde
      I_{(p)} \mid (R^T_{I,\widetilde I_{(p)}}\widetilde{\vect a}_{(p)})_j = k\}$
    \EndFor
    \State $\widetilde{\mathcal A}_{(p)} \gets \{a_k : k \in
    \widetilde I^\text{coarse}_{(p)}\}$ \Comment{Aggregate information}
    \State \Return $(\widetilde I^\text{coarse}_{(p)}, \widetilde{\mathcal A}_{(p)})$
    \EndProcedure
  \end{algorithmic}
\end{algorithm}

For each aggregate $a_k$ on process $p$, that consists of indices in
$\widetilde I_{(p)} \setminus I_{(p)}$ on the fine level, the child
node, representing that aggregate on the next coarser level, is again
associated with an index  $i \in a_k \subset \widetilde I_{(p)} \setminus I_{(p)}$. This means
that for all  
vertices in $I_{(p)}$ on the coarse level all neighbors they
depend on or influence  are also stored in process $p$.

The local prolongation operator $\widetilde P_{(p)}$ is calculated from the aggregate
information $\widetilde{\mathcal A}_{(p)}$ in accordance to
\eqref{eq:prolongation}. Let $\widetilde{\tensor A}^l_{(p)}$ be the local fine level
matrix on level $l$, then the tentative coarse level matrix is computed by the
Galerkin product $\widetilde{\tensor A}^{l+1}_{(p)}=\widetilde{\tensor
  P}_{(p)}^T\widetilde{\tensor A}^l_{(p)} \widetilde{\tensor P}_{(p)}$. To
satisfy the constraints of our local operators
\eqref{eq:tilde_A_bar_tilde}, we need to set the diagonal values to
$1$ and the off-diagonal values to $0$ for all matrix rows 
corresponding to the overlap region $\widetilde I^{l+1}_{(p)} \setminus
I_{(p)}^{l+1}$.
Due to the structure of the matrices in the hierarchy all matrix-vector
operations can be performed locally on each processor provided that
the vectors are stored consistently. 

\subsubsection{Agglomeration on Coarse Levels}
\label{sec:aggl-coarse-levels}

Note that our aggregation Algorithm
\ref{alg:parallel_coarsen} does only build aggregates within the nonoverlapping partitioning.
On the fine level, we rely on the user (or third party
software) providing our solver with a reasonable partitioning of the
global matrices and vectors onto the available processes. Often this
partitioning will not take weak connection in the matrix graph into account. 
Continuing coarsening until no further decoupled aggregation is possible, the
non-overlapping, local index sets $I_{(p)}$ are either very small
(size one in the extreme case) or the coupling between unknowns ``owned'' by
a process are all weak (which may happen e.g. in case of anisotropic problems).

In that situation there are two options: (i) either the coarsest system thus obtained is 
solved (approximately) by a single grid method (such as a preconditioned Krylov method)
or (ii) the system is redistributed to a smaller number of processors and the coarsening is
continued. The first option is viable for some problems leading to diagonally
dominant problems on the coarsest grid (e.g. certain time-dependent 
problems with time step small enough) but is in general not efficient enough. 
%The currently available supercomputers, like the Blue Gene type systems of IBM,
%already make hundreds of thousands of cores available for 
%usage. Even if the coarsest level system has only few unknowns per
%core, the global system still has several hundred
%thousand unknowns. Solving such a coarse level system in parallel
%would mean doing very few floating point operations between many
%communication steps. Therefore this computation would be limited by the
%available bandwidth and latency of the communication network. 
%To overcome these problems, there is the option to agglomerate the
%data onto fewer processes. If agglomeration is not activated, all
%processes will compute on the coarsest level and a parallel Krylov
%method preconditioned with the smoother is used as a
%solver. Otherwise, whenever the average number of unknowns per core on
%a level drops below the prescribed coarsening target a  new
%partitioning is computed using ParMETIS
In the second option a new partitioning is computed using ParMETIS
(cf. \cite{parmetisurl,parmetis}), a parallel graph partitioning
software. Our implementation supports two different choices for the input graph given to ParMETIS. 
The logically most reasonable is to use the weighted
graph of the global matrix. Its edge weights are set to $1$ for 
edges that are considered strong by our strength of connection measure and to $0$
otherwise. This tells the graph partitioning software that weak
connections can be cut at no cost and leads to partitionings that
should keep small connected regions on one process. We believe that this approach
results in sufficient coupling of strongly connected unknowns on
coarser grids. 

Unfortunately, at least as of version 3.1.1 ParMETIS uses a dense array of size $|\mathcal
P| \times |\mathcal P|$ internally to capture all possible adjacencies
between processes. This results in running out of memory on systems
with very many cores (like the IBM Blue
Gene/P). To prevent this we use the
vertex-weighted graph of the communication pattern used in the parallel
linear operator as input. Each process represents a vertex in the
graph. The weight of the vertex is the number of matrix rows stored on
this process. Edges appear only between pairs of vertices associated with
processes that exchange data. This graph is gathered on one master
process and the repartitioning is computed with the recursive graph
repartitioning routine from sequential METIS. Then the data
(matrices and vectors) of all processes associated with vertices in
one partition, is agglomerated on one process and the others become
idle on coarser levels. Obviously this is a sequential
bottleneck of our method. It will improve once massively
parallel graph partitioning tool become available.

This kind of agglomeration is repeated on subsequent levels until there is only one participating process on the
coarsest level. We can now
use a sequential sparse direct solver as coarse level solver.

\begin{figure}[htbp]
  \centering
  \includegraphics[width=.75\textwidth]{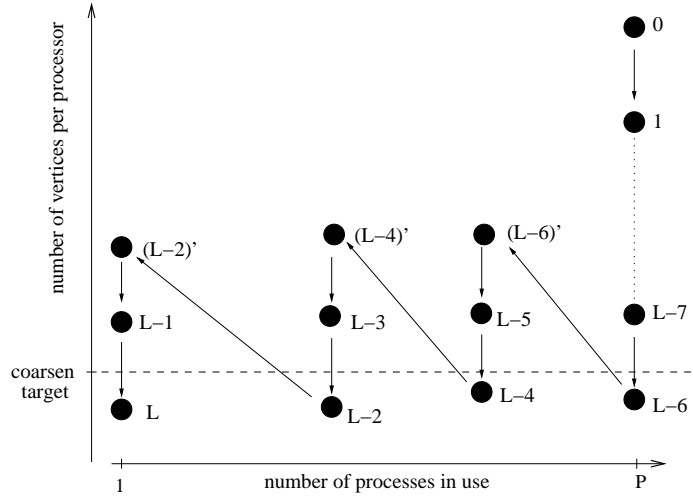}
  \caption{Data agglomeration}
  \label{fig:data_agglomeration}
\end{figure}

In Figure \ref{fig:data_agglomeration}, the interplay of the coarsening
and the data agglomeration process is sketched. Each node represents a
stored matrix. Next to it the level index is written. As before the
index $0$ denotes the finest level. Note that on
each level, where data agglomeration happens, some processes store two
matrices, an agglomerated and a non-agglomerated one. The latter is
marked with an apostrophe after the level number.

Whenever data agglomeration happened, the parallel smoothers use the not
yet agglomerated matrix. The agglomerated matrix is only needed for the coarsening to the next level. % for efficiency reasons.

\section{Implementational Details}
\label{sec:impl-deta}

The described algorithm is implemented in 
the ``Distributed and Unified Numerics Environment'' (DUNE)
\cite{DUNE,DUNEWeb,dune08-1,dune08-2}. As the components of this
library are the main cause for the good performance of our method, we
shortly introduce the two main building blocks of our AMG method: the
parallel index sets, and the ``Iterative Solver Template Library''
(ISTL) \cite{ISTL,ISTLParallel}.

\subsection{Parallel Index Sets}
\label{sec:comm-based-parall}

Our description of the parallelization in Section
\ref{sec:parallelization} is based on parallel finite index sets. This
natural representation is directly built into our software and used
for the communication. We will only shortly sketch the relevant parts
of the implementation. For a complete description of the parallel index
set software see \cite{blatt09:cplusplus_dd}.

Each process $p$ stores for each level one mapping of the
corresponding index set $\widetilde I_{(p)}$ to $\{0, \ldots,
|\widetilde I_{(p)}|-1\}$. This mapping 
allows for using the efficient local matrix and vector data structures
of ISTL to store the data and allows for direct random access. For
every entry in $\widetilde I_{(p)}$ an additional marker is stored that
lets us identify whether the index belongs to $\widetilde I_{(p)}\setminus
I_{(p)}$ or to $I_{(p)}$.  The mapping is represented by a custom
container, that provides iterators over the entries. The key type used
for this mapping is not limited to builtin integers but can be any
integral numeric type. This allows us to realize keys with enough bits
to represent integers bigger than $1.3\cdot10^{11}$. 

Using these index sets all the necessary communication patterns are
precomputed. The marker allows us for example to send data
associated with $I_{(p)}$ for every process $p\in\mathcal P$ to all
processes $q\in \mathcal P$ with 
$\widetilde I_{(q)} \cap I_{(p)} \not= \emptyset$. This kind of
communication is used for the parallel application of the linear
operator, the parallel smoother, and the communication of the
aggregate information after the decoupled aggregation of one
level. These communication patterns are implemented independent of the
communicated type. The same pattern can be used to send for example
vector entries of type double or the aggregate numbers represented by an
arbitrarily sized integer type. During the 
communication step we collect all data for each such pair $p,q$ of
processes in a buffer and send all messages simultaneously using
asynchronous communication of MPI. This keeps the number of messages
as low as possible and at the same time uses the maximal message size
possible for the problem. This reduces negative effects of network
latency. As described already in Subsections
\ref{sec:data-decomposition} and \ref{sec:smoothers} only one such
communication step is necessary for each application of the linear
operator or smoother. For the three dimensional model problems of the
next section each process sends and receives one message to and from
at most eight neighboring processes. The size of the mesage is smaller
than 152 kByte.

It is even possible to use two different distributions
$\bigcup_{p\in\mathcal P} \widetilde J_{(p)} =I$ and $\bigcup_{p\in\mathcal
  Q} \widetilde I_{(p)} =I$ as source and target of the
communication. Whenever data agglomeration occurred on a level, we have
two such distributions with $|\mathcal P|\ > |\mathcal
Q|$. To collect data we send data
associated to $I_{(p)}$ for every process $p\in\mathcal P$ to all
processes $q\in \mathcal Q$ with 
$ I_q \cap I_{(p)} \not= \emptyset$, when gathering data to
fewer processes.

\subsection{Efficient Local Linear Algebra}
\label{sec:effic-local-line}

The ``Iterative Solver
Template Library'' (ISTL) \cite{ISTL} is designed specifically for linear
systems originating from the discretization of partial differential
equations. An important application area are systems of PDEs. They
often exhibit a natural block structure. The user of our method can
choose to neglect this block structure like with most other
libraries. The linear system is then simply resembled by a sparse
matrix with scalar entries. In addition our method also supports a
block-wise treatment of the unknowns, where all unknowns associated
with the same discretization entity are grouped together. These groups
must have the same size for all entities. The couplings between the
grouped unknowns are represented by small dense matrices. The unknowns
themselves in small vectors. The size of the matrices and vectors is
known already at compile time.

ISTL offers specialized data structures for these and in addition
supports block recursive sparse matrices of arbitrary
recursion level. Using generic programming
techniques the library lets the compiler optimize the code for the
data structures during compilation. 
The available preconditioners and smoothers are implemented such that
the same code supports arbitrary block recursion levels.

Therefore our method naturally supports so-called point-based AMG,
where each matrix entry is a small dense matrix by itself. The graph
used during the coarsening in Section \ref{sec:coars-algor-based} is
the graph of the block matrix and the weight functions used in the
criterions \eqref{eq:symmetric_strong1} and \eqref{eq:symmetric_strong2} are
functions that turn the matrix blocks into scalars, such as the
row-sum or Frobenius norm. The user only has to select the appropriate
matrix and vector data structures and the smoother automatically
becomes a
block-smoother due to generic programming with templates.

\section{Numerical Results}
\label{sec:numerical-results}

In this section we present scalability results for two model
problems.  First we solve simple Laplace and Poisson problems. Then we take a
look at a heterogeneous model problem with highly variable
coefficients. We perform our test on two different hardware platforms:
a super-computer from IBM and a recent off-the-shelf Linux cluster.

The first machine is JUGENE located at the
Forschungszentrum in J\"ulich, Germany. JUGENE is a Blue Gene/P machine
manufactured by IBM that provides more than one 
petaflops as overall peak performance. Each compute node uses a 850
MHz PowerPC 450 quad-core CPU and provides 2GB of main memory with a
bandwidth of 13.6 GB/s. The main interconnect is a 3D-Torus network
for point to point message passing with a peak hardware bandwidth of
425MB/s in each direction of each torus link and a total of 5.1 GB/s
of bidirectional bandwidth for each node. Additionally there are a global
collective and a global barrier network. For comparison we also
performed some tests on helics3a at Heidelberg University, an
of-the-shelf Linux cluster consisting of 32 compute nodes with four
AMD Opteron 6212 CPUs providing eight cores, each at 2.6 GHz. Each node
utilizes 128 GB DDR3 RAM at 1333Mhz as main memory. The Infiniband
network interconnect is a Mellanox 40G QDR single port PCIe
Interconnect QSFP with 40 GB/s bidirectional bandwidths.

We start the analysis of our method by solving 
the Laplace equation, i.e. $K \equiv I$, with zero Dirichlet boundary
conditions everywhere. The results of the weak
scalability test can be found in Table \ref{tab:poisson}. For the
discretization we used a cell-centered finite volume scheme with $80^3$
cells per participating core. Note that the biggest problem computed
contains more than $1.34\cdot 10^{11}$ unknowns.

The problems are discretized
on a structured cube grid with uniform grid spacing~$h$.
We use one step of the V-cycle of the multigrid method as a preconditioner in
 a BiCGSTAB solver. For pre- and post-smoothing we apply one step of
 hybrid symmetric Gauss-Seidel. We measure the number of iterations 
(labelled It) to achieve a relative reduction of the Euclidian norm
of the residual of
$10^{-8}$. Note that BiCGSTAB naturally does apply the preconditioner two times
in  each iteration.
We measure the number of grid levels (labelled lev.), the time 
needed per iteration (labelled TIt), the time for building the AMG hierarchy 
(labelled TB), the time needed for solving the linear system
(labelled TS), and the total time needed to solution
(labelled TT) including setup and solve phase depending on the number of processors (which is proportional to the number of grid cells $1/h$ cubed). Time is always measured in seconds.

\begin{table}[htb]
  \centering
  \begin{tabular}{rr|r|rlllll}
    procs & 1/h & lev. & TB  & TS & It & TIt & TT \\\hline
    1 & 80 & 5 & 19.86  &31.91 & 8 & 3.989 & 51.77 \\
    8 & 160 & 6 & 27.7 & 46.4 & 10 & 4.64 & 74.2 \\
    64 & 320 & 7 & 74.1 &49.3 & 10 & 4.93 & 123 \\
    512 & 640 & 8 & 76.91 &60.2 & 12 & 5.017 & 137.1\\
    4096 & 1280 & 10 & 81.31 & 64.45 & 13 & 4.958 & 145.8\\
    32768 & 2560 & 11 & 92.75 & 65.55 & 13 & 5.042 & 158.3\\
    262144 & 5120 & 12 & 188.5 & 67.66 & 13 & 5.205 & 256.2 
  \end{tabular}
  \caption{Laplace Problem 3D on JUGENE: Weak Scalability}
  \label{tab:poisson}
\end{table}

Clearly, the time needed per iteration scales very well. When using
nearly the whole machine in the run with 262,144 processes, we still
reach an efficiency of about $77\%$. Due to the slight increase in
the number of iterations the efficiency of the solution phase is 
about $47\%$. Unfortunately, the hierarchy building does not scale
as well. This has different components, which can best be distinguished by an analysis of the time needed for some phases of the coarsening with agglomeration, which is
displayed in Table \ref{tab:agglo_poisson}.
%  \begin{tabular}{r|r|rrr}
%    procs & no. & TG & TM& TR \\\hline
%    1 & 1 & 0.00 & 0.00 & 0.00\\
%    8 & 1 & 0.00 & 0.39 &0.43  \\
%    64 & 1 & 0.00 & 3.56 & 3.90  \\
%    512 & 2 & 0.12 & 2.37 & 3.31\\
%    4096 & 3 & 0.71 & 3.44 & 5.80 \\
%    32768 & 3 & 7.75 & 3.31 & 15.57 \\
%    262144 & 4 & 78.25 & 5.46 & 102.73\\
%  \end{tabular}
%  \caption{Time needed for Agglomeration (Laplace 3D)}
%  \label{tab:agglo_poisson}
%\end{table}

\begin{table}[htb]
  \centering
  \begin{tabular}{r|r|r|rrr|r}
    procs & lev. &no. & TG & TM& TR \\\hline
    1 & 5 & 0 & 0.00 & 0.00 & 19.86 \\
    8 & 6 & 1 & 0.04 & 0.39 & 27.27 \\
    64 & 7 & 2 & 0.31 & 1.36 & 72.43  \\
    512 & 8 & 2 & 0.81 & 2.00 & 74.10 \\
    4096 & 10 & 3 & 2.18 & 3.04 & 75.46 \\
    32768 & 11 & 3 & 10.57 & 4.26 & 77.56 \\
    262144 & 12 & 4 & 98.23 & 4.65 & 85.71 \\
  \end{tabular}
  \caption{Time needed for Agglomeration and Coarsening(Laplace 3D)}
  \label{tab:agglo_poisson}
\end{table}
In the table the column labelled ``no.'' contains the number of data
agglomeration steps, the column labelled ``TG'' contains the time
spend in preparing the global graph on one process, partitioning it
with METIS, and creating the communication infrastructure.
The column ``TM'' contains
the time needed for redistributing the matrix data to the new
partitions, and the column labelled ``TR'' contains the total time
needed for the rest of the coarsening including the time for the
factorization of the matrix on the coarsest level using SuperLU. If we
directly agglomerate all data to one process,  we
do not use METIS as the repartitioning scheme is already known in advance. With an increasing number of processes the time
spent for computing the graph repartitioning increases much faster than the time needed
for the redistribution of the data. It turns out that this is one of the main bottlenecks, especially with very high processor numbers. However, there is also a
marked increase in the time needed for the coarsening process itself. Without 
agglomeration the creation of the coarsest matrices would take less and less time as 
the total number of entries to be aggregated decreases. However, after each redistribution step the number of matrix entries per processor still participating in the computation step is increasing again. Thus the build time has a $log(P)$ 
dependency. The time TR needed for the coarsening increases much more
whenever an additional level of agglomeration is needed.

\begin{table}[htb]
  \centering
  %% begin scaling output: weakdim=3 H=1/190problem=l correllation length=0.000976562 variance=1 mean=0 CC
\begin{tabular}{rr|rllrll}
procs & 1/h & lev. & TB & TS & It & TIt & TT \\\hline
1 & 190 & 5 & 37.97 & ~71.77 & 8 & ~8.97 & 109.7 \\
8 & 380 & 6 & 50.90 & 211.39 & 14 & 15.10 & 262.0 \\
64 & 760 & 8 & 60.00 & 243.23 & 15 & 16.20 & 303.0 \\
512 & 1520 & 9 & 66.20 & 247.75 & 15 & 16.50 & 314.0 
\end{tabular}
  \caption{Poisson Problem 3D on helics3a: Weak Scalability}
  \label{tab:poisson_helics}
\end{table}
We perform a slightly modified test on helics3a, where the Poisson problem is solved, 
described by
\begin{align*}
  -\Delta u &= (6-4\|x\|) e^{\|x\|} &&\text{in $\Omega=(0,1)^3$},\\
  u &= e^{\|x\|} &&\text{on $\partial \Omega$}
\end{align*}
As helics3a has more main memory per core it is possible to use a grid which has more than eight times as many grid cells per core than on JUGENE.

The results of the weak scalability test can be found in
Table \ref{tab:poisson_helics}. Note that despite the larger problem
per process the same number of levels in the matrix hierarchy as before is constructed.
This is equivalent to an eight times larger problem on the coarsest
level. The build time scales much better under these
circumstances. This has two reasons. First, the size of the coarse
grid problem after agglomeration is smaller compared to the large
number of unknowns per processor. Secondly, the fraction of TR needed
for the matrix factorization is higher due to the larger matrix on the
coarsest level, which reduces the influence of the graph partitioning
and redistribution. 

On helics3a with its much faster processor cores compared to JUGENE,
the memory bandwidth becomes the limiting factor for the solution phase. While for
the case of eight processes it would in principle be possible to
distribute the processes in a way that each process still has full
memory bandwidth, we did not exploit this possibility as it is very
tedious to achieve such a distribution and as it is no longer possible
for the case of 64 or more processes anyhow. With the automatic
process placement of the operating system processes will share a
memory controller already in the case of 8 processes, which is
reflected in the notable increase of the time per iteration. 
As expected the time per iteration is only slightly increasing when using
even more processes. The hierarchy
building is much less affected by the memory bandwidth limitation.  

In addition we perform a strong scalability test on helics3a where the
total problem size stays constant while the number of cores used
increases. In
this test we use decoupled coarsening until we reach the coarsening
target and then agglomerate all the data onto one process at once and
solve the 
coarse level system there. The results can be seen in Tables
\ref{tab:strong_scal} and \ref{tab:strong_eff}.  Note, that when using
512 processes our method
still has an efficiency of 27\%. Again for the time needed
for the solution phase (column TS) the biggest drop in efficiency
occurs when using eight instead of one core due to the limited memory
bandwidth.
The setup phase (column TB) is not limited as much by it and scales much better than on JUGENE.

\begin{table}[htb]
  \centering
  %% begin scaling output: strongdim=3 H=1/256problem=a correllation length=0.015625 variance=1 mean=0 CC
\begin{tabular}{r|rrrlrr}
procs & TB & TS & It & TIt & TT \\\hline
1 & 102.60 & 166.20 & 6 & 27.70 & 268.80 \\
8 & 14.00 & 35.80 & 8 & 4.47 & 49.80 \\
64  & 2.02 & 5.06 & 9 & 0.56 & 7.08 \\
512 & 0.73 & 1.16 & 8 & 0.15 & 1.89 
\end{tabular}
%% end scaling output
  \caption{Poisson Problem 3D on helics3a: Strong Scalability}
  \label{tab:strong_scal}
\end{table}

\begin{table}[htb]
  \centering
%% begin scaling output: strongdim=3 H=1/256problem=a correllation length=0.015625 variance=1 mean=0 CC
\begin{tabular}{r|llll}
procs & TB & TS & TIt & TT \\\hline
8 & 0.92 & 0.58 & 0.77 & 0.68\\
64 & 0.79 & 0.51 & 0.77 & 0.59\\
512 & 0.27 & 0.28 & 0.37 & 0.28
\end{tabular}
%% end scaling output
  \caption{Poisson Problem 3D on helics3a: Strong Efficiency}
  \label{tab:strong_eff}
\end{table}

The last model problem we investigate is the diffusion problem
$$
\nabla \cdot (k(\vect x)\nabla x) = f
$$
on the unit cube $[0,1]^3$ with Dirichlet boundary conditions and
jumps in the diffusion coefficient as  proposed in
\cite{Griebel.Metsch.Schweitzer:2007}. Inside the unit cube, a smaller
cube with width $.8$ is centered such that all faces are parallel
to the faces of the enclosing cube. The diffusion coefficient $k$ in this
smaller cube is
$10^3$. Outside of the small cube $k=1$ holds except for cubes
with width $0.1$ that are placed in the corners of the unit cube. There
the diffusion coefficient is $10^{-2}$. Again we use a cell-centered discretization
scheme and the same settings as before for the AMG. The results
of a weak scalability test on 64 racks of JUGENE can be found in Table
\ref{tab:griebel}.
\begin{table}[htb]
  \centering
  \begin{tabular}{rr|r|rlllll}
    procs & 1/h & lev. & TB & TS & It & TIt & TT \\\hline
1 & 80 & 5 & 19.88 & 36.27 & 9 & 4.029 & 56.15\\
8 & 160 & 6 & 27.8 & 48.9 & 10 & 4.89 & 76.7\\
64 & 320 & 7 & 74.4 & 59.6 & 12 & 4.96 & 134 \\
512 & 640 & 8 & 78.04 & 72.67 & 14 & 5.191 & 150.7 \\
4096 & 1280 & 10 & 89.72 & 73.37 & 14 & 5.241 & 163.1 \\
32768 & 2560 & 11 & 94.48 & 104.2 & 20 & 5.21 & 198.7 \\
262144 & 5120 & 12 & 186.2 & 85.87 & 16 & 5.367 & 272.1
 \end{tabular}
  \caption{Heterogeneous Diffusion Problem 3D on JUGENE: Weak Scalability}
  \label{tab:griebel}
\end{table}

Compared to the Poisson problem the number of iterations increases
more steeply due to the jumps in the diffusion coefficient. Again this is
not due to the parallelization but due to the nature of the
problem. The time needed for building the matrix hierarchy as well as
the time needed for one iteration scale as for the previous
problems. Compared to the AMG method used in
\cite{Griebel.Metsch.Schweitzer:2007} on an (now outdated) Blue Gene/L
our method scales much better when used on Blue Gene/P. In small parts this
might due to the new architecture and the bigger problem size per core
used. But this cannot explain all the difference in the
scaling behavior. Additionally, the large problem size is only possible
because of the smaller memory foot-print of our method. 

\section{Summary and Conclusion}
\label{sec:summary}
 
We have presented a  parallel algebraic multigrid algorithm
based on non-smoothed aggregation. During the setup phase it uses an
elaborate heuristic aggregation algorithm to account for highly variable
coefficients that appear in many application areas. Due to its
simple piecewise constant interpolation between the levels, the memory
consumption of the method is rather low and allows for solving
problems with more than $10^{11}$ unknowns using 64 racks of an IBM
Blue Gene/P.
The parallelization of the solution phase scales well for up to nearly
300,000 cores. Although there is a sequential bottleneck in the setup
phase of the method due to the lack of scalable parallel graph
partitioning software, the method still scales very well in terms of
total time to solution. For comparison see \cite{baker12:amg_exascale}
where during a weak scalability test for the Laplace problem on an IBM
Blue Gene/P the total solution time for interpolation AMG increases by
more than a factor 2 when going from 128 to 128,000 processes. In
contrast, for our method with the same increase in total solution time
we can go from 64 to up to 262,144 processes during weak scaling.
Additional comparisons made by M\"uller and Scheichl \cite{Eike2013} support
the claim that our implementation scales at least as well as the Ruge-St\"{u}ben
type AMG presented in \cite{baker12:amg_exascale} also for anisotropic problems
and that floating-point performance is comparable.

We also have shown that our solver scales reasonably well even for hard
problems that have highly variable coefficients. Even for modern
clusters consisting out of multicore machines the method scales very
well and is only limited by the available memory bandwidth per core.

Once scalable parallel graph partition software is available, the
bottleneck of the sequential graph partitioning will disappear
rendering the method even more scalable.
\bibliography{paamg,literatur}
\bibliographystyle{siam}
\end{document}